\documentclass{amsart}
\usepackage{amssymb}
\usepackage{amsfonts}
\usepackage{latexsym}

\newtheorem{theorem}{Theorem}[section]
\newtheorem{lemma}[theorem]{Lemma}

\newtheorem{corollary}[theorem]{Corollary}
\theoremstyle{definition}
\newtheorem{definition}[theorem]{Definition}

\newtheorem{remark}[theorem]{Remark}

\newcommand{\Fun}{\text{Fun}}
\newcommand{\Irr}{\text{\rm Irr}}
\newcommand{\FPdim}{\text{\rm FPdim}}

\newcommand{\Hom}{\text{Hom}}

\newcommand{\Rep}{\text{Rep}}

\newcommand{\ot}{\otimes}

\newcommand{\ben}{\begin{enumerate}}
\newcommand{\een}{\end{enumerate}}

\newcommand{\Vect}{\text{Vec}}

\newcommand{\cC}{{\mathcal C}}

\newcommand{\cM}{{\mathcal M}}
\newcommand{\be}{{\bf 1}}

\newcommand{\BZ}{{\mathbb Z}}

\newcommand{\BC}{{\mathbb C}}

\newcommand{\iHom}{\underline{\mbox{Hom}}}

\begin{document}

\title{Classification of Fusion Categories of
Dimension $pq$}
\author{Pavel Etingof}
\address{Department of Mathematics, Massachusetts Institute of Technology,
Cambridge, MA 02139, USA} \email{etingof@math.mit.edu}

\author{Shlomo Gelaki}
\address{Department of Mathematics, Technion-Israel Institute of
Technology, Haifa 32000, Israel}
\email{gelaki@math.technion.ac.il}

\author{Viktor Ostrik}
\address{Department of Mathematics, Massachusetts Institute of Technology,
Cambridge, MA 02139, USA} \email{ostrik@math.mit.edu}

\maketitle

\section{Introduction}
A fusion category over $\Bbb C$
is a $\Bbb C$-linear semisimple rigid tensor category with finitely many
simple objects and finite dimensional spaces of morphisms,
such that the neutral object is simple (see \cite{eno}).
To every fusion category, one can attach a positive number,
called the Frobenius-Perron (FP) dimension of this category
(\cite{eno}, Section 8). It is an interesting and challenging
problem to classify fusion categories of a given
FP dimension $D$. This problem is easier if $D$ is an integer,
and for integer $D$ its complexity increases with the
number of prime factors in $D$. Specifically, fusion categories of
FP dimension $p$ or $p^2$
where $p$ is a prime were classified in \cite{eno},
Section 8. The next level of complexity is
fusion categories of FP dimension $pq$, where
$p<q$ are distinct primes. In this case, the classification has
been known only in the case when the category admits a fiber
functor, i.e. is a representation category of a Hopf algebra
(\cite{gw,eg1}).

In this paper we provide a complete classification
of fusion categories of FP dimension $pq$,
thus giving a categorical generalization of \cite{eg1}. As a
corollary we also obtain the classification of semisimple
quasi-Hopf algebras of dimension $pq$.
A concise formulation of our main result is:

\begin{theorem}
Let $\mathcal{C}$ be a fusion category over $\Bbb C$
of FP dimension $pq$, where $p<q$ are distinct primes. Then
either $p=2$ and $\mathcal{C}$ is a Tambara-Yamagami category of
dimension $2q$ (\cite{ty}), or $\mathcal{C}$
is group-theoretical in the sense of \cite{eno}.
\end{theorem}

The organization of the paper is as follows.

In Section 2 we recall from \cite{eno} and \cite{o} some
facts about fusion categories that will be used below.

In Section 3 we classify fusion categories $\mathcal{C}$
of dimension $pq$ which contain objects of non-integer
dimensions. They exist only for $p=2$ and are the
Tambara-Yamagami categories \cite{ty}; there are four such
categories for each $q$.

In Section 4 we classify fusion categories
of dimension $pq$ where all simple objects are invertible.
This reduces to computing the cohomology groups
$H^3(G,\Bbb C^*)$, where $G$ is a group of order $pq$.

In Section 5 we deal with the remaining (most difficult) case,
when $\mathcal{C}$ has integer dimensions of simple objects,
but these dimensions are not all equal to $1$.
We show, by generalizing the methods of \cite{eg1},
 that in this case $q-1$ is divisible by $p$, and
the simple objects of $\mathcal{C}$ are
$p$ invertible objects and $\frac{q-1}{p}$ objects
of dimension $p$.

In Section 6, we classify fusion categories ${\mathcal C}$
whose simple objects are
$p$ invertible objects and $\frac{q-1}{p}$ objects
of dimension $p$. Namely, we show that they are group-theoretical
in the sense of \cite{eno}, which easily yields their
full classification.

As a by-product, the method of Section 6 yields a classification
of finite dimensional semisimple quasi-Hopf (in particular, Hopf)
algebras whose irreducible representations have dimensions $1$ and
$n$, such that the 1-dimensional representations form a cyclic
group of order $n$. All such quasi-Hopf algebras turn out to be
group-theoretical. This is proved in Section 7. We also classify
fusion categories whose invertible objects form a cyclic group of
order $n>1$ and which have only one non-invertible object of
dimension $n$.

We note that all constructions in this paper are done over the field
of complex numbers.

{\bf Acknowledgments.} The first author was partially supported by
the NSF grant DMS-9988796. The first author partially conducted
his research for the Clay Mathematics Institute as a Clay
Mathematics Institute Prize Fellow. The second author's research
was supported by Technion V.P.R. Fund - Dent Charitable Trust- Non
Military Research Fund, and by THE ISRAEL SCIENCE FOUNDATION
(grant No. 70/02-1). The third author's work was partially
supported by the NSF grant DMS-0098830.

\section{Preliminaries}

Let $\mathcal{C}$ be a fusion category over $\mathbb{C}$ and
let $K(\mathcal{C})$ denote its Grothendieck ring. Let
$\Irr(\mathcal{C})$ be the (finite) set of isomorphism classes of
simple objects in $\mathcal{C}$ and let $V\in \Irr(\mathcal{C})$.

\begin{definition}\cite{eno}
(i) The Frobenius-Perron (FP) dimension of $V$, $\FPdim(V)$, is the
largest positive eigenvalue of the matrix of multiplication by $V$ in
$K(\mathcal{C})$.

(ii) The FP dimension of $\mathcal{C}$ is
$\FPdim(\mathcal{C})=\sum_{V\in \Irr(\mathcal{C})} \FPdim(V)^2$.

\end{definition}

Let $Z(\mathcal{C})$ be the Drinfeld center of $\mathcal{C}$. Then by
Proposition 8.12 in \cite{eno},
$\FPdim(Z(\mathcal{C}))=\FPdim(\mathcal{C})^2$.

If $\mathcal{C}$ is a full tensor subcategory of $\mathcal{D}$
then $\FPdim(\mathcal{D})/\FPdim(\mathcal{C})$ is an algebraic
integer, so in particular, if the two dimensions are integers then
the ratio is an integer (\cite{eno}, Proposition 8.15).
Moreover, if $\mathcal{M}$ is a full module subcategory
of $\mathcal{D}$ over $\mathcal{C}$, then
the same is true about $\FPdim(\mathcal{M})/\FPdim(\mathcal{C})$,
where $\FPdim(\mathcal{M})$ is the sum of squares
of the Frobenius-Perron dimensions of simple objects of $\mathcal{M}$
(\cite{eno}, Remark 8.17).

By Proposition 8.20 in \cite{eno}, if $\mathcal{C}=\bigoplus_{g\in G}
\mathcal{C}_g$ is faithfully graded by a finite group $G$ then
$\FPdim(\mathcal{C}_g)$ are equal for all $g\in G$ and $|G|$
divides $\FPdim(\mathcal{C})$.

There is another notion of a dimension for fusion categories
$\mathcal{C}$; namely, their global dimension $\dim(\mathcal{C})$
(see \cite{eno}, Section 2). The global dimension of $\mathcal{C}$
may be different from its FP dimension. However, categories
$\mathcal{C}$ for which these two dimensions coincide are of great
interest. They are called {\em pseudounitary} \cite{eno}. One of
the main properties of pseudounitary categories is that they admit
a unique pivotal structure. This extra structure allows one to define
categorical dimensions of simple objects of $\mathcal{C}$, which by
Proposition 8.23 in \cite{eno}, coincide with their FP dimensions.
One instance in which it is guaranteed that $\mathcal{C}$ is
pseudounitary is when $\FPdim(\mathcal{C})$ is an integer
(Proposition 8.24 in \cite{eno}). If $\mathcal{C}$ has integer FP dimension,
then the dimensions of simple objects in $\mathcal{C}$ are integers
and square roots of integers (\cite{eno}, Proposition 8.27).

An important special case of fusion categories
with integer FP dimension are categories
in which the FP dimensions of all simple objects are integers.
It is well known (see e.g. Theorem 8.33 in \cite{eno}) that
this happens if and only if the category is equivalent to
$\Rep(H)$ where $H$ is a finite dimensional semisimple quasi-Hopf
algebra.

By Corollary 8.30 in \cite{eno}, if $\FPdim(\mathcal{C})$ is equal
to a prime number $p$, then $\mathcal{C}$ is equivalent
$\Rep({\rm Fun}(\mathbb{Z}/p\mathbb{Z}))$
with associativity defined by a cocycle $\xi$ representing $\omega\in
H^3(\mathbb{Z}/p\mathbb{Z},\mathbb{C}^*)=\mathbb{Z}/p\mathbb{Z}$.

An important class of fusion categories with integer FP dimensions
of simple objects is the class of {\em group theoretical fusion
categories}, introduced and studied in \cite{eno}, \cite{o}. These
are categories associated with quadruples $(G,B,\xi,\psi)$,
where $G$ is a finite group, $B$ is a subgroup of $G$, $\xi\in
Z^3(G,\mathbb{C}^*)$, and $\psi\in C^2(B,\mathbb{C}^*)$ such that
$\xi_{|B}=d\psi$. Namely, let $\Vect_{G,\xi}$ be the
category of finite dimensional $G-$graded vector spaces with
associativity defined by $\xi$. Let $\Vect _{G,\xi}(B)$ be
the subcategory of $\Vect _{G,\xi}$ of objects graded by $B$.
Consider the twisted group algebra $A:=\mathbb{C}^{\psi}[B]$. It
is an associative algebra in $\Vect _{G,\xi}(B)$, since
$\xi_{|B}=d\psi$. Then $\mathcal{C}(G,B,\xi,\psi)$ is
defined to be the category of $A-$bimodules in $\Vect
_{G,\xi}$. Such a category is called {\em group-theoretical}.

Note that the data $(\xi,\psi)$ is not uniquely determined by the category.
Namely, there are two transformations of $(\xi,\psi)$ which leave the
category unchanged:

1) $\xi\to \xi+d\phi,\psi\to \psi+\phi|_B$, $\phi\in
C^2(G,\Bbb C^*)$; and

2) $\xi\to \xi$, $\psi\to \psi+d\eta$, $\eta\in C^1(B,\Bbb
C^*)$.

Thus the essential data is the cohomology class of $\xi$ (which
must vanish when restricted to $B$) and an element $\psi$ of a
principal homogeneous space (torsor) $T_\xi$ over the group
${\rm Coker}(H^2(G,\Bbb C^*)\to H^2(B,\Bbb C^*))$ (\cite{eno},
Remark 8.39).

Proposition 8.42 in \cite{eno} gives a simple
characterization of group-theoretical
fusion categories. Namely, a fusion
category $\mathcal{C}$ is group-theoretical if and only if it is
dual to a pointed category (= category whose all simple objects are
invertible) with respect to some indecomposable module category over
$\mathcal{C}$. For the definitions of all unfamiliar terms we
refer the reader to \cite{eno}, \cite{o}.

\section{Categories with non-integer dimensions}

 Throughout the paper we
will consider a fusion category
$\mathcal{C}$ over $\mathbb{C}$ of
FP dimension $pq$, where $p$ and $q$ are
distinct primes, such that $p<q$. As explained in Section 2,
such a category is pseudo-unitary,
and hence admits a canonical pivotal structure, in which
categorical dimensions of objects coincide with their
FP dimensions. Thus from now on we will refer to FP dimensions
simply as ``dimensions''.

By Proposition 8.27
in \cite{eno}, the dimensions of simple objects in $\mathcal{C}$
may be integers or square roots of integers.

\begin{theorem}\label{notall}
If $\mathcal{C}$ contains a simple object whose dimension is not
an integer then $p=2$
and $\mathcal{C}$ is equivalent to a Tambara-Yamagami category
of dimension $2q$ \cite{ty}.
\end{theorem}

\begin{proof} Let $\mathcal{C}_{ad}$ be the full tensor subcategory
of $\mathcal{C}$ generated by the constituents in $X\ot X^* \in
\mathcal{C}$, for all simple $X\in \mathcal{C}$. By
Proposition 8.27 in \cite{eno}, the dimensions of objects in
$\mathcal{C}_{ad}$ are integers. Since $\mathcal{C}$ has objects
whose dimension is not an integer, $\mathcal{C}\ne
\mathcal{C}_{ad}$. Thus by Proposition 8.15 in \cite{eno}, the
dimension $d$ of $\mathcal{C}_{ad}$ is an integer dividing $pq$ and
less than $pq$, so it is either $1$ or $p$ or $q$. If $d=1$
then all objects of $\mathcal{C}$ are invertible, contradiction.
Thus $d=p$ or $q$.

By Proposition 8.30 in
\cite{eno}, $\mathcal{C}_{ad}$ is $\Rep(\Fun(\mathbb{Z}/d\mathbb{Z}))$
with associativity defined by a 3-cocycle (i.e. the simple objects are
$\chi^i$, $i=0,...,d-1$).

Let $\mathcal{C}'$ be the full subcategory
of $\mathcal{C}$ consisting of objects with integer dimension.
We claim that it is a tensor subcategory of $\mathcal{C}$.
Indeed, if $X\oplus Y$ has integer dimension, then so do
$X$ and $Y$ (sum of square roots of positive integers is an
integer only if each summand is an integer).

Thus, $\mathcal{C}'$ is a proper tensor subcategory of $\mathcal{C}$
containing $\mathcal{C}_{ad}$.
So by Proposition 8.15 in \cite{eno}, $\mathcal{C'}$ coincides with
$\mathcal{C}_{ad}$.

Now, let $L\in \mathcal{C}$ be a simple object, such that $L\notin
\mathcal{C}_{ad}$. Then one has $L\ot L^*=\oplus_i \chi^i$.
Indeed, since $\dim(L)>1$, $L\otimes L^*$ must contain simple
objects other than $\bold 1$ (which it contains with multiplicity
$1$). But $L\otimes L^*\in {\mathcal C}_{ad}$, so the other
constituents could only be $\chi^i$, and they can only occur with
multiplicity $1$.

Thus, the dimension of $L$ is $\sqrt{d}$, and $\chi^i\ot L=L\ot
\chi^i=L$.

Let us now show that $L$ is unique. Let $L,M$ be two such
simple objects (in $\mathcal{C}$ but not in $\mathcal{C}_{ad}$).
Then $M\ot L^*$ has dimension $d$, so it lies in
$\mathcal{C}_{ad}$. Hence, $\chi^i$ occurs in $M\otimes L^*$ for
some $i$. So $\chi^i\ot L=M$, and hence $L=M$.

Thus, $p=2$, $d=q$, and
the dimension of $\mathcal{C}$ is $2q$. Moreover, the simple
objects of $\mathcal{C}$
are $q$ invertible objects and one object of dimension
$\sqrt{q}$. Hence $\mathcal{C}$
is a Tambara-Yamagami category, as desired.
\end{proof}

\section{Categories with 1-dimensional simple objects}

It is well known \cite{k} that fusion categories
with 1-dimensional simple objects are classified by pairs
$(G,\omega)$, where $G$ is a finite group, and $\omega\in H^3(G,\Bbb C^*)$.
Namely, the category ${\mathcal C}(G,\omega)$ attached to such a pair is the
category of representations of the function algebra ${\rm Fun}(G)$ with
associativity defined by a cocycle $\xi$ representing $\omega$.
Thus, in order to classify such categories of dimension $pq$,
it is sufficient to classify pairs $(G,\omega)$
with $|G|=pq$.

There are two kinds of groups of order $pq$: the cyclic group
$G=\Bbb Z/pq\Bbb Z$, and the nontrivial semidirect product $G
=\BZ/p\BZ \ltimes \BZ/q\BZ$ (which exists and is unique if and
only if $q-1$ is divisible by $p$). In the first case, it is well
known that $H^3(G,\Bbb C^*)=\Bbb Z/pq\Bbb Z$. So it remains to
consider the second case, $G=\BZ/p\BZ \ltimes \BZ/q\BZ$, where
$q-1$ is divisible by $p$, and the action of $\BZ/p\BZ$ on
$\BZ/q\BZ$ is nontrivial.

\begin{lemma}\label{hs}
For $i>0$ we have $H^i(G, \BC^*)=H^i(\BZ/p\BZ,\BC^*)\oplus
(H^i(\BZ/q\BZ,\BC^*))^{\BZ/p\BZ}$. \end{lemma}

\begin{proof} The Hochschild-Serre spectral sequence has second term
\linebreak $E_2^{ij}= H^i(\BZ/p\BZ,H^j(\BZ/q\BZ,\BC^*))$. For any
$j>0$, $H^j(\BZ/q\BZ,\BC^*)$ is a $q-$group and thus $E_2^{ij}=0$
if both $i,j$ are nonzero. Moreover, all differentials are zero
since they would map $q-$groups to $p-$groups. Thus the spectral
sequence collapses and the lemma is proved.
\end{proof}

\begin{corollary}\label{coho} One has
$H^3(G,\BC^*)=H^3(\BZ/p\BZ,\BC^*)=\BZ/p\BZ$ if $p\ne 2$ and
\linebreak $H^3(G,\BC^*)=H^3(\BZ/2\BZ,\BC^*)\oplus
H^3(\BZ/q\BZ,\BC^*)= \BZ/2\BZ \oplus \BZ/q\BZ$ if $p=2$.
\end{corollary}

\begin{proof} By Lemma \ref{hs},
it is enough to check that the action of $\BZ/p\BZ$
on $H^3(\BZ/q\BZ,\BC^*)$ is nontrivial for odd $p$ and trivial for $p=2$.
For this, observe that $ H^3(\BZ/q\BZ,\BC^*)=H^4(\BZ/q\BZ,\BZ)=
(H^2(\BZ/q\BZ,\BZ))^{\otimes 2}=(H^1(\BZ/q\BZ,\BC^*))^{\otimes 2}=
(\Hom(\BZ/q\BZ,\BC^*))^{\otimes 2}$
as $\Bbb Z/p\Bbb Z$-modules, and the claim is proved.
\end{proof}

\section{Categories with integer dimensions, not all equal to $1$}

In this section we will prove the following result.

\begin{theorem}\label{1andp}
Let ${\mathcal C}$ be a fusion category of dimension $pq$
with integer dimensions
of simple objects, not all equal to $1$.
Then $q-1$ is divisible by $p$, and
the simple objects of ${\mathcal C}$ are $p$ invertible objects
and $\frac{q-1}{p}$ objects of dimension $p$.
\end{theorem}

The proof of this theorem occupies the rest of the section.

In the proof we will assume that $pq\ne 6$, because in the case
$pq=6$ the theorem is trivial.

By Theorem 8.33 in
\cite{eno}, $\mathcal{C}$ is equivalent to $\Rep(H)$, where $H$ is
a finite dimensional semisimple quasi-Hopf algebra. Let
$Z(\mathcal{C})$ be the Drinfeld center of $\mathcal{C}$, then
$Z(\mathcal{C})$ is equivalent to $\Rep(D(H))$, where $D(H)$ is
the double of $H$ \cite{hn}.

\begin{lemma}\label{1pq}
The simple objects in ${\rm Rep}(D(H))$ have dimension $1$, $p$ or
$q$.
\end{lemma}

\begin{proof}
Since $\mathcal{C}$ is pivotal, $Z(\mathcal{C})$ is modular, and
the result follows from Lemma 1.2 in \cite{eg2} (see also
Proposition 3.3 in \cite{eno}).
\end{proof}

\begin{lemma}\label{nt1}
$D(H)$ admits nontrivial 1-dimensional representations;
i.e., the group of grouplike elements $G(D(H)^*)$
of the coalgebra $D(H)^*$ is nontrivial.
\end{lemma}

\begin{proof}
Assume the contrary. Let $m$
be the number of $p$-dimensional representations of $D(H)$, and $n$
the number of its $q$-dimensional representations.
Then by Lemma \ref{1pq}, one has
$1+mp^2+nq^2=p^2q^2$, which implies that $m>0$.

Let $V$ be a $p-$dimensional irreducible representation of $D(H).$
Then $V\ot V^*$ is a
direct sum of the trivial representation $\Bbb C,$
$a$ $p-$dimensional irreducible
representations of $D(H)$
and $b$ $q-$dimensional irreducible representations of $D(H).$ Therefore, we
have: $p^2=1+ap+bq.$ Clearly $b>0.$ Let $W$ be a $q-$dimensional
irreducible
representation of $D(H)$ such that $W\subset V\ot V^*.$ Since
$$
0\ne
{\rm Hom}_{D(H)}(V\ot
V^*,W)={\rm Hom}_{D(H)}(V,W\ot V),
$$ we have that $V\subset W\ot V.$ Since
$W\ot V$ has no 1-dimensional
constituent (because $\dim V\ne \dim W$), $\dim (W\ot V)=pq$
and $W\ot V$ contains a
$p-$dimensional irreducible representation of $D(H)$,
from dimension counting it follows that
$W\ot V=V_1\oplus\cdots \oplus V_q$ where $V_i$ is a $p-$dimensional
irreducible representation of $D(H)$ with $V_1=V.$

We wish to show that for any $i=1,\dots,q$, $V_i=V.$ Suppose on
the contrary that this is not true for some $i.$ Then $V\ot V_i^*$
has no 1-dimensional constituent, hence it must be a direct sum of
$p$ $p-$dimensional irreducible representations of $D(H).$
Therefore, $W\ot (V\ot V_i^*)$ has no 1-dimensional constituent.
But, $(W\ot V)\ot V_i^*=(V_i\oplus\cdots )\ot V_i^*=\Bbb
C\oplus\cdots$, which is a contradiction.

Therefore, $W\ot V=qV.$
Hence ${\rm Hom}_{D(H)}(V\ot
V^*,W)$ is a q-dimensional space, i.e. $p^2=\dim(V\otimes V^*)\ge q^2$,
contradiction.
\end{proof}

\begin{lemma}\label{inj} (i) The natural map $r: G(D(H)^*)\to G(H^*)$ is
injective.

(ii) $|G(D(H)^*)|=p$ or $q$, and thus $r$ is an isomorphism.
\end{lemma}

\begin{proof} Assume the contrary. Then there is a non-trivial cyclic
subgroup $L$ in $G(D(H)^*)$ which maps trivially to $G(H^*)$. This
means that the category $\mathcal{C}$ is faithfully
$L^\vee$-graded, by Proposition 5.10 in \cite{eno}. So,
$\mathcal{C}$ is a direct sum of $\mathcal{C}_\gamma$, $\gamma\in
L^\vee$, and the dimension of $\mathcal{C}_\gamma$ is $s:=pq/|L|$,
which is $1$ or $p$ or $q$. If $s=1$, all simple objects of
$\mathcal{C}$ are invertible, which is a contradiction. If $s=p$
or $s=q$ then $\mathcal{C}_0$ is a fusion category of prime
dimension. So by Proposition 8.30 in \cite{eno}, ${\mathcal C_0}=
{\mathcal C}(\Bbb Z/s\Bbb Z,\omega)$ and the $\mathcal{C}_\gamma$'s
are module categories over it. If $\mathcal{C}_\gamma$ has a
non-1-dimensional object $V$ then $\chi\ot V=V\ot \chi=V$ for
$\chi\in \mathcal{C}_0$ (otherwise, $\dim(\mathcal{C}_\gamma)$ will
be greater than $s$), so $\dim(V)$ is divisible by $s$ (Remark
8.17 in \cite{eno}), i.e. $\dim(\mathcal{C}_\gamma)\ge s^2$,
contradiction. Thus, again, all simple objects of $\mathcal{C}$
are invertible, which is a contradiction.

Thus $r$ is injective. By Proposition
8.15 in \cite{eno} and Lemma \ref{nt1},
this means that $|G(D(H)^*)|$ is either $p$ or $q$, as desired.
\end{proof}

\begin{lemma}\label{nq}
$|G(D(H)^*)|=p$, and $q-1$ is divisible by $p$.
\end{lemma}

\begin{proof}
Set $m:=|G(D(H)^*)|$. Let $\mathcal{D}\subset\mathcal{C}$ be the
subcategory generated by the invertible objects in $\mathcal{C}$
(it is of dimension $m$). Then we have $\mathcal{C}\subset
\mathcal{C}\boxtimes \mathcal{D}^{op}\subset \mathcal{C}\boxtimes
\mathcal{C}^{op}$. Taking the dual of this sequence with respect
to the module category $\mathcal{C}$ (see \cite{eno}, Sections 5
and 8), we get a sequence of surjective functors
$Z(\mathcal{C})\to \mathcal{E}\to \mathcal{C}$, where the
dimension of $\mathcal{E}$ is $mpq$. We can think of the category
$\mathcal{E}$ as representations of a quasi-Hopf subalgebra
$B\subset D(H)$, containing $H$, of dimension $mpq$.

Let $\chi$ be a 1-dimensional representation of $D(H)$. Let
$J(\chi):=D(H)\ot_{B}\chi$ be the induced module. By Schauenburg's
freeness theorem \cite{S1} (see also \cite{eno}, Corollary 8.9), it
has dimension $pq/m$.

For any two 1-dimensional representations $\chi,\chi'$ of $D(H)$,
we have $$\Hom_{D(H)}(J(\chi),\chi')=\Hom_{B}(\chi,\chi').$$ By
Lemma \ref{inj} (i), this is zero if $\chi\ne \chi'$ and $\Bbb C$
if $\chi=\chi'$. Thus the only 1-dimensional constituent of
$J(\chi)$ as a $B$-module is $\chi$, and it occurs with
multiplicity $1$.

Assume that $m=q$. Then the dimension of $J(\chi)$ is
$p$. Since other
constituents of $J(\chi)$ can only have dimensions $1,p,q$, and
$p<q$, we get that $J(\chi)$ is a sum of characters $\chi'$ of
$D(H)$, which is a contradiction. So,
by Lemma \ref{inj} (ii), $m=p$, as desired.

In this case, $\dim (J(\chi))=q$, and $J(\chi)$ must be $\chi$
plus sum of $p-$dimensional simple modules, whose number is then
$(q-1)/p$. Thus, $q-1$ is divisible by $p$, as desired.
\end{proof}

\begin{lemma}\label{start}
Let $V,U$ be $p$-dimensional representations of $D(H)$.
If $V\otimes U$ contains a 1-dimensional representation $\chi$ of $D(H)$,
then it contains another 1-dimensional representation.
\end{lemma}

\begin{proof}
Without loss of generality we can assume that $\chi$ is
trivial and
hence that $U=V^*.$ Otherwise we can replace $U$ with $U\ot \chi^{-1}.$

Suppose on the contrary that $V\ot V^*$ does not contain a non-trivial
$1-$dimensional representation. Then $V\ot V^*$ is a
direct sum of the trivial representation $\Bbb C,$ $p-$dimensional irreducible
representations of $D(H)$
and $q-$dimensional irreducible representations of $D(H).$ Therefore we
have that $p^2=1+ap+bq.$ Since by Lemma \ref{nq}, $q-1$ is
divisible by $p$, we find that $b+1$ is divisible by $p$. So
$b\ge p-1$, and hence $p^2\ge 1+(p-1)q$,
i.e. $p+1\ge q$. Thus $p=2, q=3$, and $pq=6$.
But we assumed that it is not the case, so we have a contradiction.
\end{proof}

\begin{lemma}\label{np}
The algebra $D(H)$ has $p^2-p$ $q$-dimensional irreducible
representations and $(q^2-1)/p$ $p$-dimensional irreducible
representations. Moreover, the direct sums of 1-dimensional and
p-dimensional irreducible representations form a tensor
subcategory ${\mathcal F}$ in ${\rm Rep}(D(H))$ of dimension
$pq^2$.
\end{lemma}

\begin{proof}
Let $a$ and $b$ be the numbers of $p$-dimensional and $q$-dimensional
irreducible representations of $D(H)$. Then by Lemma \ref{nq} (i),
$ap^2+bq^2=p^2q^2-p$. This equation clearly has a unique
nonnegative integer solution $(a,b)$. By Lemma \ref{nq} (ii),
this solution is $a=(q^2-1)/p, b=p^2-p$.

Let us now prove the second statement.
We wish to show that if $V$ and $U$ are two irreducible representations
of $D(H)$ of dimension $p$ then $V\ot U$ is a direct sum of $1-$dimensional
irreducible representations of $D(H)$ and $p-$dimensional irreducible
representations of $D(H)$ only. Indeed,
by Lemma \ref{start}, either $V\ot U$ does not contain any $1-$dimensional
representation or it must contain at least two different $1-$dimensional
representations. But if it contains two different $1-$dimensional
representations, then (since the $1-$dimensional representations of $D(H)$
form a cyclic group of order $p$) $V\ot U$ contains all
the $p$ $1-$dimensional
representations of $D(H).$ We conclude that either $p^2=mp+nq$ or
$p^2=p+mp+nq.$ At any rate $n=0,$ and the result follows.
\end{proof}

Therefore, the subcategory ${\rm Rep}(D(H))$ generated by the $1$
and $p-$dimensional irreducible representations of $D(H)$ is the
representation category of a quotient quasi-Hopf algebra $A$ of
$D(H)$ of dimension $pq^2$.

\begin{lemma}\label{sur}
The composition map $H\to D(H)\to A$ is injective.
Thus $H$ is a quasi-Hopf subalgebra in $A$.
\end{lemma}

\begin{proof} Assume that the composition map is not injective.
Then the image of this map is a nontrivial quotient of $H$. The
image definitely contains the subalgebra $A_0$ in $A$
corresponding to the invertible objects. This quasi-Hopf
subalgebra is $p$-dimensional, while $H$ is $pq$-dimensional, so
by Schauenburg's theorem \cite{S2}, see also \cite{eno},
Proposition 8.15, the image
must coincide with $A_0$. On the other hand,
by Schauenburg freeness theorem \cite{S1}, $D(H)$ is a free
left $H$-module of rank $pq$. Since the projection
$D(H)\to A$ is a morphism of left $H$-modules,
we find that $A$ is generated by $pq$ elements as a left
$A_0$-module. Hence, the dimension of $A$ is at most $p^2q$.
On the other hand, we know that this dimension is $pq^2$, a
contradiction.
\end{proof}

\begin{lemma}\label{irre}
Any irreducible representation $V$ of $H$ which is not 1-dimensional
has dimension $p$.
\end{lemma}

This lemma clearly completes the proof of the theorem.

\begin{proof}
It is clear from Lemma \ref{sur} that this dimension is at most
$p$ (as any simple $H$-module occurs as a constituent in a simple
$A$-module). On the other hand, we claim that $V$ is stable under
tensoring with 1-dimensional representations. Indeed, assume not,
and let $W$ be an irreducible representation of $A$ whose
restriction to $H$ contains $V$. Since $W$ is $p$-dimensional, and
contains all $\chi^j\otimes V$ (where $\chi$ is a non-trivial
1-dimensional representation of $H$), by dimension counting
we get that $p\ge p\dim
(V)$, a contradiction.

But now by Remark 8.17 in \cite{eno}, the dimension of $V$ is
divisible by $p$. We are done.
\end{proof}

\section{Categories with integer dimensions}

In this section we will prove that
any fusion category of dimension $pq$
with integer dimensions of objects is group
theoretical, and will classify such categories.
Before we do so, we need to prove two lemmas.

\begin{lemma} \label{fp}
Let $\cC$ be a fusion category and let $A\in \cC$ be an
indecomposable semisimple algebra. Then for any right $A-$module $M$ and
left $A-$module $N$ one has ${\rm FPdim}(M\otimes_AN)=\frac{{\rm FPdim}(M)
{\rm FPdim}(N)}{{\rm FPdim}(A)}$. \end{lemma}

\begin{proof} Let $M_i, i\in I$ be the collection of simple right $A-$modules
and let $N_j, j\in J$ be the collection of simple left
$A-$modules. It is clear that it is enough to prove the lemma for
$M=M_i$ and $N=N_j$. Note that the vector ${\rm FPdim}(M_i)$
(resp. ${\rm FPdim}(N_j)$) is the Frobenius-Perron eigenvector
(see \cite{eno}) for the module category of right $A-$modules
(resp. left $A-$modules). For any left $A-$module $N$ the vector
${\rm FPdim}(M_i\otimes_AN)$ is also the Frobenius-Perron
eigenvector and thus is proportional to ${\rm FPdim}(M_i)$.
Similarly, for any right $A-$module $M$ the vector ${\rm
FPdim}(M\otimes_AN_j)$ is proportional to ${\rm FPdim}(N_j)$. Thus
${\rm FPdim}(M_i\otimes_AN_j)=\alpha {\rm FPdim}(M_i){\rm
FPdim}(N_j)$ for some constant $\alpha$. Finally, by choosing
$M=N=A$ we find out that $\alpha =1/{\rm FPdim}(A)$. The lemma is
proved.
\end{proof}

Let $n>1$ be an integer. Let $\cC$ be a fusion category and $\chi
\in \cC$ be a nontrivial invertible object such that $\chi^{n}=
\be$. Assume that $\cC$ contains a simple object $V$ such that
$\chi \ot V= V\ot \chi = V$. This implies that $A:=\be \oplus \chi
\oplus \chi^{2}\oplus \ldots \oplus \chi^{n-1}$ has a unique
structure of a semisimple algebra in $\cC$. Indeed, the existence
of $V$ implies that there is a fiber functor (= module category
with one simple object $V$) on the category generated by
$\{\chi^i\}$, i.e. the 3-cocycle of this category is trivial and
thus it is the representation category of the Hopf algebra ${\rm
Fun}(\mathbb{Z}/n\mathbb{Z})$. Then the dual to this Hopf algebra
is the algebra $A$.

Assume additionally that for any simple object $X$ of $\cC$ we
have either ${\rm FPdim}(X)=n$ or $X$ is isomorphic to $\chi^{i}$
for some $i$.

\begin{lemma}
\label{inv} Let $M$ be a simple $A-$bimodule such that ${\rm
Hom}_\cC (V,M)\ne 0$. Then $M= V$ as an object of $\cC$. In
particular $M$ is invertible in the tensor category
of $A-$bimodules.
\end{lemma}

\begin{proof}
Assume first that $M$ is a simple right $A-$module such that
$\Hom_\cC (V,M)\ne 0$. Then $\Hom_A(V\otimes A,M)=\Hom_\cC
(V,M)\ne 0$ and hence $M$ is a direct summand of $V\otimes A=
V^{\oplus n}$. On the other hand it is obvious that the object $V$
has $n$ different structures of an $A-$module. Thus we have proved
that any simple right (and similarly left) $A-$module $M$ with
$\Hom_\cC (V,M)\ne 0$ is isomorphic to $V$ as an object of $\cC$.

Now let $M$ be a simple right $A-$module such that $M= V$ as an
object of $\cC$. Let us calculate $\iHom(M,M)$. By example 3.19 in
\cite{eo}, $\iHom(M,M)=(M\otimes_A{}^*M)^*$ and by Lemma \ref{fp},
${\rm FPdim}(\iHom(M,M))=n$. Clearly, $\be \subset \iHom(M,M)$ and
thus $\iHom(M,M)=\be \oplus \chi \oplus \ldots \oplus \chi^{n-1}$.
In particular $\chi \otimes M= M$ as right $A-$modules. Choose
such an isomorphism; after normalizing we can consider it as a
structure of a {\em left} $A-$module on $M$ commuting with the
structure of a right $A-$module. In other words, $M$ has $n$
different structures of an $A-$bimodule (for a fixed right
$A-$module structure). Thus altogether we constructed $n^2$
different structures of an $A-$bimodule on the object $V$.
Finally, any simple $A-$bimodule $M$ with $\Hom_\cC (V,M)\ne 0$ is
a direct summand of $A\otimes V\otimes A = V^{\oplus n^2}$. The
lemma is proved.
\end{proof}

Now we are ready to state and prove the main result of this section.

\begin{theorem}\label{integer}
Let $p<q$ be primes. Then
any fusion category ${\mathcal C}$ of Frobenius-Perron dimension $pq$
with integer dimensions of
simple objects is group-theoretical. More specifically, it is
equivalent to one in the following list:

(i) A category with 1-dimensional simple objects (these are described in
Section 4).

(ii) ${\rm Rep}(G)$ where $G =\BZ/p\BZ \ltimes \BZ/q\BZ$ is a
non-abelian group.

(iii)
If $p=2$, the
category ${\mathcal C}(G,\Bbb Z/2\Bbb Z,\xi,\psi)$
(see \cite{eno}, Section 8.8; \cite{o})
where $G$ is
the nonabelian group
$\BZ/2\BZ \ltimes \BZ/q\BZ$, and
$\xi\in Z^3(G,\Bbb C^*)=\Bbb Z/2q\Bbb Z$ is a
cocycle which represents a cohomology class of order $q$,
and $\psi$ is determined by $\xi$.
In this case ${\mathcal C}$ is not a representation category of a Hopf
algebra.
\end{theorem}

\begin{proof}
Assume that not all simple objects of ${\mathcal C}$ are
invertible. To prove the first statement, observe that by Theorem
\ref{1andp}, the assumptions of Lemma \ref{inv} are satisfied.
Hence Lemma \ref{inv} applies. Thus, any simple $A$-bimodule $M$
containing a $p$-dimensional representation $V$ is invertible. On
the other hand, it is clear that any simple $A$-bimodule which
involves only $\chi^{i}$ must be isomorphic to $A$. Thus, any
simple $A$-bimodule is invertible. In other words, the dual
category ${\mathcal C}_{\Rep(A)}^*=A-{\rm bimod}$ has only
invertible simple objects. So ${\mathcal C}$ is group-theoretical,
as desired.

Let us now prove the second statement. The category ${\mathcal
C}_{\Rep(A)}^*$ is of the form ${\mathcal C}(G,\omega)$, where $G$
is a group of order $pq$, and $\omega\in H^3(G,\Bbb C^*)$. So
${\mathcal C}$ is of the form ${\mathcal C}(G,B,\xi,\psi)$,
where $B$ is a subgroup of $G$, $\xi$ a 3-cocycle representing $\omega$
and $\psi$ is a 2-cochain on $B$.

It is easy to check that if the category ${\mathcal
C}(G,B,\xi,\psi)$ has non-1-dimensional simple objects then $G$ is
the nonabelian group $\Bbb Z/p\Bbb Z\ltimes \Bbb Z/q\Bbb Z$, and
$B=G$ or $B=\Bbb Z/p\Bbb Z$.

Further, the cocycle $\xi$ must be trivial on $B$, and $\psi$
is determined by $\xi$ up to equivalence (see Lemma \ref{hs}).
If $B=G$ then this implies that we can set $\xi=1$, and we are
in case (ii). Suppose that $B=\Bbb Z/p\Bbb Z$. If $p$ is odd, then
Lemma \ref{coho} implies that we can set $\xi=1$. In this case
$\mathcal{C}$ is the representation category of the Kac algebra
attached to the exact factorization $G=(\Bbb Z/p\Bbb Z)(\Bbb
Z/q\Bbb Z)$. It is easy to see that this Kac algebra is isomorphic
to the group algebra of $G$ as a Hopf algebra, so we are still in
case (ii).

If $p=2$ and $\xi=1$ in cohomology,
we are in case (ii) as well, for the
same reason. If $p=2$ and $\xi\ne 1$ in cohomology,
then we are in case (iii).
It follows from \cite{o} that in this case $\mathcal{C}$ does not admit
fiber functors.
We are done.
\end{proof}

\begin{remark}
Let $\xi$ be a 3-cocycle on $\mathbb{Z}/q\mathbb{Z}$ representing
a nontrivial cohomology class. Since $\mathbb{Z}/2\mathbb{Z}$ acts
trivially on $H^3(\mathbb{Z}/q\mathbb{Z},\mathbb{C}^*)$,
and since $2$ is relatively prime to $q$,
$\xi$ can be chosen to be invariant under $\mathbb{Z}/2\mathbb{Z}$.
Let $\Phi$ be an associator in $\Fun(\mathbb{Z}/q\mathbb{Z})^{\ot
3}$ corresponding to $\xi$. Then $(\mathbb{C}[\mathbb{Z}/2\mathbb{Z}]
\ltimes \Fun(\mathbb{Z}/q\mathbb{Z}),\Phi)$, with the usual coproduct,
is a finite dimensional semisimple quasi-Hopf algebra $H$.
Then $\Rep(H)$ is
a category from case (iii), and any category of case (iii)
(there are two of them up to equivalence) is obtained in this way.
\end{remark}

\begin{remark}
Theorem \ref{integer} implies in particular the classification of
semisimple quasi-Hopf algebras of dimension $pq$, where $p$ and
$q$ are distinct primes.
\end{remark}

\begin{remark} In the case $pq=6$, Theorem \ref{integer} was
proved by T.~Chmutova. Namely, she discovered
that besides categories whose simple objects are invertible,
there are exactly three 6-dimensional categories with integer
dimensions of simple objects:
the category of representations of $S_3$ (case (ii)) and two
additional categories with the same Grothendieck ring
(case (iii)).
\end{remark}

\section{Categories with simple objects
of dimension $1$ and $n$}

Let $n>1$ be an integer. Let $N$ be a finite group with a
fixed-point-free action of $\BZ/n\BZ$ and let $\omega \in H^3(N,\BC^*)$
be an invariant class under the $\BZ/n\BZ-$action. Since $n$ and
$|N|$ are coprime there exists a 3-cocycle $\xi$ representing
$\omega$ and invariant under the $\BZ/n\BZ-$action. Let $\Phi$ be
an associator in $\Fun(N)^{\ot 3}$ corresponding to $\xi$. Then
$(\mathbb{C}[\mathbb{Z}/n\mathbb{Z}] \ltimes \Fun(N),\Phi)$, with
the usual coproduct, is a finite dimensional semisimple quasi-Hopf
algebra $H$. It is easy to
see that any simple $H-$module has dimension 1 or $n$. The
following theorem gives an abstract characterization of quasi-Hopf
algebras constructed in such a way.

\begin{theorem}\label{new} Let $\cC$ be a fusion category such that

(i) Invertible objects of $\cC$ form a cyclic group of order $n$.

(ii) For any simple object $X\in \cC$ either ${\rm FPdim}(X)=1$ or
${\rm FPdim}(X)=n$, and $\cC$ contains at least one simple object
of FP dimension $n$.

Then there exists a finite group $N\ne \lbrace{1\rbrace}$
with a fixed-point-free
action of $\BZ/n\BZ$ and a $\BZ/n\BZ-$invariant class $\omega \in
H^3(N,\BC^*)$ such that $\cC$ is equivalent to ${\rm Rep}(H)$
where $H$ is the quasi-Hopf algebra constructed above.
\end{theorem}

\begin{proof} Let $m$ be the number of simple objects $X$ in $\cC$
with ${\rm FPdim}(X)=n$. Then ${\rm FPdim}(\cC)=n(mn+1)$.

Let $V$ be an $n$-dimensional simple object
of $\mathcal{C}$. Then $V\ot V^*$ contains
the neutral object, so by dimension counting it must contain all
$1-$dimensional objects. Thus $V$ is stable under tensoring with
$1-$dimensional objects. Hence, we are in the conditions of Lemma
\ref{inv}. Thus we see that the
category $\Rep(A)$ has exactly $mn+1$ simple objects (namely, the regular
module and $n$ structures of an $A$-module on each simple
noninvertible object of $\cC$), and the dual category
$\cC^*_{\Rep(A)}$ has only invertible objects. Thus
$\cC^*_{\Rep(A)}=\cC(G,\omega)$ for some finite group $G$ and
$\omega \in H^3(G,\BC^*)$. Moreover, it is clear from the
classification of module categories over $\cC(G,\omega)$ that
$\Rep(A)$, as a module category over $\cC^*_{\Rep(A)}$, is of the
form $\cM (G,B,\xi,\psi)$ where $B\subset G$ is a cyclic
subgroup of order $n$ such that $\xi|_B$ is trivial, see
\cite{o}. Thus, $\psi$ is determined by $\xi$ up to equivalence
(as $H^2(\Bbb Z/n\Bbb Z, \Bbb C^*)=0$), and $\cC =
(\cC^*_{\Rep(A)})^*_{\Rep(A)}=\cC(G,B,\xi,\psi)$. The simple
objects in the category $\cC(G,B,\xi,\psi)$ are classified by
pairs $(g,\lambda)$ where $g\in G$ and $\lambda$ is an irreducible
representation of $B\cap gBg^{-1}$ (see \cite{o}), and the
Frobenius-Perron dimension of the simple object corresponding to a
pair $(g,\lambda)$ is $|B:B\cap gBg^{-1}|\dim (\lambda)$. Thus
conditions (i), (ii) translate to the following: $B\cap
gBg^{-1}=1$ for any $g\not \in B$. In other words, $G$ is a
Frobenius group (see e.g. \cite{g}). Thus $G\simeq B\ltimes N$ for
some normal subgroup $N\subset G$ and the action of $B$ on $N$ is
fixed-point-free (see {\em loc. cit.}). Furthermore,
$H^3(G,\BC^*)=H^3(B,\BC^*)\oplus H^3(N,\BC^*)^B$ by the same
argument as in Lemma \ref{hs}. Clearly the subgroup of $\omega \in
H^3(G,\BC^*)$ such that $\omega|_B=1$ is identified with
$H^3(N,\BC^*)^B$. Thus, by the Frobenius theorem, $G=B\ltimes N$.
The Theorem is proved.
\end{proof}

\begin{corollary} Let $H$ be a semisimple Hopf algebra
with $1$-dimensional and $n$-dimensional irreducible
representations, such that $G(H^*)$ is a cyclic group of order
$n$. Then $H=\Bbb C[B]\ltimes {\rm Fun}(N)$ is the Kac algebra
attached to the exact factorization $G=BN$, where $B=\Bbb Z/n\Bbb
Z$, $N$ is a group with a fixed-point-free action of $B$, and
$G=B\ltimes N$.
\end{corollary}

\begin{proof} In the Hopf algebra case
the category $\cC$ admits a fiber functor, i.e. a module category
with only one simple object. Thus, by \cite{o}, there exists a
subgroup $P$ of $G$ such that $G=BP$, $\omega|_P=1$. Clearly, $P$
contains $N$, so $\omega=1$, $\psi=1$, and we are done.
\end{proof}

\begin{remark} Recall that the
famous Thompson's Theorem states that the group $N$ above is
nilpotent (see e.g. \cite{g}).
\end{remark}

Now consider the special case of Theorem \ref{new} when $m=1$. Let
$X$ denote the non-invertible object and $\chi^i, i=0, \ldots,
n-1,$ denote the invertible objects of $\cC$. Then obviously the
multiplication in the category $\cC$ is given by
$$\chi \otimes X=X\otimes \chi =X,\; X\otimes X=(n-1)X\oplus \chi^0\oplus
\ldots \oplus \chi^{n-1}.$$

\begin{corollary} Let $\cC$ be a fusion category such that the invertible
objects of $\cC$ form a cyclic group of order $n>1$ and $\cC$ has
only one non-invertible object of dimension $n$. Then $n+1=p^a$ is
a prime power. If $n=2$ there are three such categories, if $n=3$
or $7$ there are two such categories, and for all other
$n=p^a-1>1$ there is exactly one such category -- the category of
representations of the semi-direct product ${\mathbb F}_{p^a}^*
\ltimes {\mathbb F}_{p^a}$.
\end{corollary}

\begin{remark} We thank R. Guralnick for help in the proof of the
corollary.
\end{remark}

\begin{proof} In this case the group $N$ above is of order $n+1$ and the
group $B$ acts simply transitively on the non-identity elements of
$N$. Thus all non-identity elements of $N$ have the same order and
hence $N$ is a $p-$group. Consequently, any element of $N$ is
conjugated to some central element and hence $N$ is abelian.
Henceforth, $N$ is an elementary abelian group of order $q=p^a$.
The cyclic group $B$ acts irreducibly on $N$, hence by Schur's
Lemma, $N$ is identified with a one dimensional vector space over
the finite field ${\mathbb F}_q$, and $B$ is identified with
$GL_1({\mathbb F}_q)={\mathbb F}_q^*$.

The following statement is well known:

\begin{lemma} Let $V$ be an elementary abelian $p-$group. Consider
$H^i(V,\BC^*)$ as a functor in the variable $V$. Then we have

(i) $H^1(V,\BC^*)=V^*$.

(ii) $H^2(V, \BC^*)=\wedge^2V^*$.

(iii) There is an exact sequence of $GL(V)$-modules
$0\to S^2V^*\to H^3(V,\BC^*)\to \wedge^3V^*\to 0$.
\end{lemma}

Here $S^\bullet V^*$ is the symmetric algebra of the space $V^*$;
that is, the algebra generated by $v\in V^*$ subject to the
relations $v_1v_2=v_2v_1$ for any $v_1, v_2\in V^*$. Similarly,
$\wedge^\bullet V^*$ is generated by $v\in V^*$ subject to the
relations $v^2=0$.

\begin{proof} Items (i) and (ii) are well known. We prove (iii). Recall that
$H^3(V,\BC^*)=H^4(V,\BZ)$. It follows from the Kunneth formula
that $H^{>0}(V,\BZ)$ is annihilated by the multiplication by $p$.
Thus an exact sequence $0\to \BZ \to \BZ \to \BZ/p\BZ \to 0$ for
any $i\ge 1$ gives an exact sequence $0\to H^i(V,\BZ)\to
H^i(V,\BZ/p\BZ)\to H^{i+1}(V,\BZ)\to 0$. It is well known (see
e.g. \cite{Be}) that
$$H^\bullet (V,\BZ/p\BZ)=\left\{ \begin{array}{cc}S^{2\bullet}V^*\otimes
\wedge^\bullet V^*& \mbox{for } p>2\\ S^\bullet V^*& \mbox{for }
p=2\end{array} \right.$$ Thus for $p>2$ we have
$H^3(V,\BZ/p\BZ)=V^*\otimes V^*\oplus \wedge^3V^*$ and for $p=2$,
$H^3(V,\BZ/p\BZ)=S^3V^*$. For $p>2$ one observes that the image of
$\wedge^2V^*=H^3(V,\BZ)$ lies inside $V^*\otimes V^*$, since the
scalar matrices act by different characters on $\wedge^2V^*$ and
$\wedge^3V^*$. Also, $V^*\otimes V^*/\wedge^2V^*=S^2V^*$ and we
are done. For $p=2$ we get that there is an embedding
$\wedge^2V^*\subset S^3V^*$ and $H^4(V,\BZ)=S^3V^*/ \wedge^2V^*$.
Consider the obvious surjection $S^3V^*\to \wedge^3V^*$. Since
$\wedge^2V^*$, $\wedge^3V^*$ are simple non-isomorphic
$GL(V)-$modules, the submodule $\wedge^2V^*$ is in the kernel of
this surjection. On the other hand, it is easy to see that the
kernel is identified with $V^*\otimes V^*$ via the map $x\otimes
y\mapsto x^2y$. Finally, one observes that $V^*\otimes V^*$ has a
unique copy of the simple module $\wedge^2V^*$
spanned by tensors of the form
$x\otimes y+y\otimes x$, and $V^*\otimes V^*/\wedge^2V^*=S^2V^*$.

The lemma is proved.
\end{proof}

Now, one deduces easily that in our situation $H^3(N,\BC^*)^B$ is
nontrivial if and only if $q=3,4,8$. Indeed, let $\alpha$ be a
generator of ${\mathbb F}_q^*$. Then the operator of
multiplication by $\alpha$ in the vector space $V^*={\mathbb F}_q$
has eigenvalues $\alpha, Fr(\alpha)=\alpha^p, \ldots,
Fr^{a-1}(\alpha)=\alpha^{p^{a-1}}$. The eigenvalues for the action
on $S^2V^*$ (resp. $\wedge^3V^*$) are $\alpha^{p^i+p^j},\; 0\le
i\le j\le a-1$ (resp. $\alpha^{p^i+p^j+p^k},\; 0\le i<j<k\le
a-1$). Thus we have an eigenvalue 1 on $S^2V^*$ (resp.
$\wedge^3V^*$) if and only if $p^i+p^j=p^a-1$ for some $0\le i\le
j\le a-1$ (resp. $p^i+p^j+p^k=p^a-1$ for $0\le i<j<k\le a-1$), and
the statement follows. In all cases the space $H^3(N,\BC^*)^B$ is
one dimensional over the prime field and the case $q=3$ was
already considered in Theorem \ref{integer}. Finally, note that
the category $\cC :=\Rep({\mathbb F}_q^*\ltimes {\mathbb F}_q)$
satisfies the conditions of the corollary. This completes the
proof of the corollary.
\end{proof}


\begin{thebibliography}{AEG}

\bibitem[B]{Be} D. Benson, Representations and Cohomology,
{\em Cambridge Studies in Advanced Mathematics}, 1991.

\bibitem[EG1]{eg1} P. Etingof and S. Gelaki, Semisimple Hopf algebras of
dimension $pq$ are trivial, {\em Journal of Algebra} {\bf 210}
(1998), 664--669.

\bibitem[EG2]{eg2} P. Etingof and S. Gelaki, Some properties of
finite-dimensional semisimple Hopf algebras, {\em Mathematical
Research Letters} {\bf 5} (1998), 191--197.

\bibitem[ENO]{eno} P. Etingof, D. Nikshych and V. Ostrik, On
fusion categories, {\em preprint}, math.QA/0203060.

\bibitem[EO]{eo} P. Etingof and V. Ostrik, Finite tensor
categories, {\em preprint}, math.QA/0301027.

\bibitem[G]{g} D. Gorenstein, Finite Groups, Harper \& Row, Publishers,
New York--London 1968.

\bibitem[GW]{gw} S. Gelaki and S. Westreich, On semisimple Hopf
algebras of dimension $pq,$ {\em Proceedings of the AMS}, {\bf
128} (2000), no.1, 39--47.

\bibitem[HN]{hn} F. Hausser and F. Nill, Doubles of quasi-quantum
groups, {\em Comm. Math. Phys.} {\bf 199} (1999), no. 3, 547--589.

\bibitem[K]{k} C. Kassel, Quantum Groups, Springer, New York, 1995.

\bibitem[O]{o} V. Ostrik, Boundary conditions for holomorphic
orbifolds, math.QA/0202130.

\bibitem[S1]{S1} P. Schauenburg,
A quasi-Hopf algebra freeness theorem, math.QA/0204141.

\bibitem[S2]{S2} P. Schauenburg, Quotients of finite quasi-Hopf
algebras, math.QA/0204337.

\bibitem[TY]{ty} D. Tambara and S. Yamagami,
Tensor categories with fusion rules of self-duality for finite
abelian groups, {\em J. Algebra} {\bf 209} (1998), no. 2,
692--707.

\end{thebibliography}
\end{document}